\documentclass[12pt]{article}
\usepackage{amsmath,amssymb,amsbsy,amsfonts,amsthm,latexsym,
               amsopn,amstext,amsxtra,euscript,amscd}
\newtheorem{theorem}{Theorem}[section]
\newtheorem{lemma}[theorem]{Lemma}

\newtheorem{corollary}[theorem]{Corollary}
\newtheorem{conjecture}[theorem]{Conjecture}
\theoremstyle{definition}

\def\cB{\mathcal B}

\def\cL{\mathcal L}
\def\cM{\mathcal M}
\def\cN{\mathcal N}

\def\cV{\mathcal V}
\def\cW{\mathcal W}

\def\R{{\mathbb R}}

\def\eps{\varepsilon}
\def\e{\mathbf{e}}

\def\mand{\qquad \mbox{and} \qquad}

\def\\{\cr}
\def\({\left(}
\def\){\right)}
\def\[{\left[}
\def\]{\right]}
\def\<{\langle}
\def\>{\rangle}
\def\fl#1{\left\lfloor#1\right\rfloor}
\def\rf#1{\left\lceil#1\right\rceil}
\def\le{\leqslant}
\def\ge{\geqslant}

\begin{document}

\title{\sc Density of non-residues in Burgess-type
intervals and applications}

\author{
{\sc W.~D.~Banks} \\
{Department of Mathematics} \\
{University of Missouri} \\
{Columbia, MO 65211 USA} \\
{\tt bbanks@math.missouri.edu} \\
\and
{\sc M.~Z.~Garaev} \\
{Instituto de Matem{\'a}ticas}\\
{ Universidad Nacional Aut\'onoma de M{\'e}xico} \\
{C.P. 58089, Morelia, Michoac{\'a}n, M{\'e}xico} \\
{\tt garaev@matmor.unam.mx} \\
\and
{\sc D.~R.~Heath-Brown} \\
{Mathematical Institute}\\
{24--29, St. Giles'} \\
{Oxford, OX1 3LB} \\
{\tt rhb@maths.ox.ac.uk} \\
\and
{\sc I.~E.~Shparlinski} \\
{Department of Computing}\\
{Macquarie University} \\
{Sydney, NSW 2109, Australia} \\
{\tt igor@ics.mq.edu.au}}

\date{\empty}

\maketitle

\newpage

\begin{abstract}
We show that for any fixed $\eps>0$, there are numbers $\delta>0$
and $p_0\ge 2$ with the following property: for every prime $p\ge
p_0$ and every integer $N$ such that $p^{1/(4\sqrt{e}\,)+\eps}\le
N\le p$, the sequence $1,2,\ldots,N$ contains at least $\delta N$
quadratic non-residues modulo~$p$. We use this result to obtain
strong upper bounds on the sizes of the least quadratic non-residues
in Beatty and Piatetski--Shapiro sequences.
\end{abstract}


\section{Introduction}

In 1994 Heath-Brown conjectured the existence of an absolute
constant~$c>0$ such that, for all positive integers $N$ and all
prime numbers $p$, the interval $[1,N]$ contains at least $cN$
quadratic residues modulo~$p$. This conjecture has been established
by Hall~\cite{Hall}. In the seminal work of Granville and
Soundararajan~\cite{GrSo} it has been shown that if $N$ is
sufficiently large, then for every prime $p$ more than $17.15\%$ of
the integers in $[1,N]$ are quadratic residues modulo~$p$. On the
other hand, for any fixed positive integer~$N$ there exist
infinitely many primes $p$ such that the interval $[1,N]$ is free of
quadratic non-residues modulo~$p$; see~\cite{GrRing} for a more
precise statement. In particular, complete analogues of the results
of Hall~\cite{Hall} and of Granville and Soundararajan~\cite{GrSo}
are not possible in the case of quadratic non-residues.

In the present paper we show that for any given $\eps>0$ there
exists a constant $c(\eps)>0$ with the following property: for every
sufficiently large prime $p$ and every integer $N$ in the range
$p^{1/(4\sqrt{e}\,)+\eps}\le N\le p$, the interval $[1,N]$ contains
at least $c(\eps) N$ quadratic non-residues modulo~$p$. This is the
partial analogue of Hall's result for quadratic non-residues in
\emph{Burgess-type} intervals. We recall that the celebrated result
of Burgess~\cite{Burg} states that the least positive quadratic
non-residue modulo~$p$ is of size $O\(p^{1/(4\sqrt{e}\,)+\eps}\)$
for any given $\eps>0$, and the constant $1/(4\sqrt{e}\,)$ has never
been improved.

We apply our result on the density of non-residues to obtain strong
upper bounds on the sizes of the least quadratic non-residues in
Beatty and Piatetski-Shapiro sequences, which substantially improve
all previously known results for these questions.

\section{Statement of results}
\label{sec:notate}

For an odd prime $p$, we use $(\cdot|p)$ to denote the Legendre
symbol modulo~$p$, and we put
$$
S_p(x)=\sum_{n\le x}(n|p)\qquad(x\ge 1).
$$

\begin{theorem}
\label{thm:PosPropNonres} For every $\eps>0$ there exists $\delta>0$
such that, for all sufficiently large primes $p$, the bound
$$
|S_p(N)|\le(1-\delta)N
$$
holds for all integers $N$ in the range $p^{1/(4\sqrt{e}\,)+\eps}\le
N\le p$.
\end{theorem}

For two fixed real numbers $\alpha$ and $\beta$, the corresponding
\emph{non-homogeneous Beatty sequence} is the sequence of integers
defined by
$$
\cB_{\alpha,\beta}=\(\fl{\alpha n+\beta}\)_{n=1}^\infty.
$$
Beatty sequences appear in a variety of apparently unrelated
mathematical settings, and because of their versatility, the
arithmetic properties of these sequences have been extensively
explored in the literature; see, for example,
\cite{Abe,Beg,Kom1,Kom2,OB,Tijd} and the references contained
therein.

For each prime $p$, let $N_{\alpha, \beta}(p)$ denote the least
positive integer $n$ such that $\fl{\alpha n + \beta}$ is a
quadratic non-residue modulo~$p$ (we formally put $N_{\alpha,
\beta}(p)=\infty$ if no such integer exists). Below, we show that
Theorem~\ref{thm:PosPropNonres} can be applied to establish the
following Burgess-type bound, which substantially improves earlier
results in~\cite{BaSh1,BaSh2,Gar,Preo1,Preo2,Preo3}:

\begin{theorem}
\label{thm:nonres} Let $\alpha,\beta$ be fixed real numbers with
$\alpha$ irrational. Then, for every $\eps>0$ the bound
$$
N_{\alpha, \beta}(p) \le p^{1/(4\sqrt e\,)+\eps}
$$
holds for all sufficiently large primes $p$.
\end{theorem}

We remark that the irrationality of $\alpha$ is essential to our
argument. Even in the ``simple'' case $\alpha=3,\,\beta=1$, we have
not been able to improve upon the inequality
$$
N_{3,1}(p)\le p^{1/4+o(1)}
$$
which follows from the Burgess bound on the relevant character sum.

Next, let $N_c(p)$ be the least positive integer $n$ such that
$\fl{n^c}$ is a quadratic non-residue modulo~$p$. It is easy to show
that $N_c(p)$ exists for any non-integer $c>1$. For values of $c$
close to $1$, good upper bounds for $N_c(p)$ have been obtained
in~\cite{Gar,LW}. Here, we establish a much stronger bound by
appealing to Theorem~\ref{thm:PosPropNonres}. It is formulated in
terms of \emph{exponent pairs}, we refer to~\cite{GrKol,Huxley,Ivic,
Rob1, Rob2, Sarg} for their exact definition  and properties.

\begin{theorem}
\label{thm:nonresPiat} Let $(\kappa,\lambda)$ be an exponent pair,
and suppose that
$$
1<c<1+\frac{1-\lambda}{2\kappa-\lambda+3}.
$$
Then, for every $\eps>0$ the bound
$$
N_{c}(p) \le p^{1/(4(2-c)\sqrt e\,)+\eps}
$$
holds for all sufficiently large primes $p$.
\end{theorem}

The classical exponent pair $(\kappa,\lambda) = (1/2,1/2)$ implies
that Theorem~\ref{thm:nonresPiat} is valid for $c$ in the range
$1<c<8/7$. Graham's optimization algorithm (see~\cite{Gr, GrKol})
extends this range to
$$
1<c<1+\frac{1-R}{2-R}=1.14601346\cdots,
$$
where $R=0.8290213568\cdots$ is Rankin's constant. Note that as
$c\to 1^+$ our upper bound for $N_c(p)$ tends to the Burgess bound,
which illustrates the strength of our estimate.

\section{Proofs}

\subsection{Proof of Theorem~\ref{thm:PosPropNonres}}
\label{sec:proof}

We can assume that $0<\eps\le 0.01$. In view of the identities
$$
\#\{n\le x~:~(n|p)=\pm 1\}=\sum_{n\le
x}\tfrac12(1\pm(n|p))=\tfrac12(\fl{x}\pm S_p(x))\qquad(x\ge 1),
$$
and taking into account the result of Hall~\cite{Hall} mentioned
earlier, it suffices to establish only the lower bound
$$
\#\{n\le N~:~(n|p)=-1\}\ge \tfrac12\,\delta\,N
$$
with $N$ in the stated range.

By the character sum estimate of Hildebrand~\cite{Hild} (which
extends the range of validity of the Burgess bound~\cite{Burg}) it
follows that $S_p(p^{1/4})=o(p^{1/4})$ as $p\to\infty$; therefore,
$$
\#\{n\le p^{1/4}~:~(n|p)=-1\}=(0.5+o(1))p^{1/4}.
$$
Since every non-residue $n$ is divisible by a prime non-residue $q$,
we have
$$
(0.5+o(1))p^{1/4}\le \sum_{n\le
p^{1/4}}\sum_{\substack{q\,\mid\,n\\
(q|p)=-1}}1\le\sum_{\substack{q\le p^{1/4}\\
(q|p)=-1}}\frac{p^{1/4}}{q},
$$
and thus
$$
0.5+o(1)\le\sum_{j=1}^s \frac{1}{q_j}
+\sum_{p^{1/(4\sqrt{e}\,)+0.5\eps}<q\le p^{1/4}}\frac{1}{q},
$$
where $q_1<\cdots<q_s$ are the prime quadratic non-residues
modulo~$p$ that do not exceed $p^{1/(4\sqrt{e}\,)+0.5\eps}$. Using
Mertens' formula (see~\cite[Theorem~427]{HarWr}), we bound the
latter sum by
$$
\sum_{p^{1/(4\sqrt{e}\,)+0.5\eps}<q\le
p^{1/4}}\frac{1}{q}=\log\(\frac{\log p^{1/4}}{\log
p^{1/(4\sqrt{e}\,)+0.5\eps}}\)+O\(\frac{1}{\log p}\)\le
0.5-2\,\varepsilon,
$$
where the inequality holds for all sufficiently large $p$.
Consequently,
$$
\sum_{j=1}^s \frac{1}{q_j}\ge\varepsilon
$$
if the prime $p$ is large enough.

For each $j=1,\ldots,k$, let $\cN_j$ denote the set of positive
quadratic residues modulo~$p$ which do not exceed $N/q_j$. From the
result of Granville and Soundararajan~\cite{GrSo} we have
$$
\#\cN_j\ge\frac{0.1N}{q_j}\qquad (j=1,\ldots,s).
$$
In particular, if $q_1\le\eps^{-1}$, then the numbers
$$
\{q_1n~:~n\in\cN_1\}
$$
are all positive non-residues of size at most $N$, and the theorem
follows from the lower bound $\#\cN_1\ge 0.1\eps N$.

Now suppose that $q_1>\eps^{-1}$. In this case, we can choose $k$
such that
$$
\eps\le\sum_{\ell=1}^{k}\frac{1}{q_\ell} \le 2\,\eps.
$$
For each $j=1,\ldots,s$, let $\cM_j$ be the set of numbers in
$\cN_j$ that are not divisible by any of the primes $q_1,\ldots,
q_k$; then
$$
\#\cM_j\ge\#\cN_j-\sum_{\ell=1}^k\frac{N}{q_jq_\ell}\ge
\frac{(0.1-2\,\eps)N}{q_j}\ge\frac{0.09N}{q_j},
$$
where we have used the fact that $\eps\le 0.01$ for the last
inequality. It is easy to see that the numbers of the form $q_jn$
with $j\in \{1,\ldots,k\}$ and $n\in \cM_j$ are distinct
non-residues of size at most $N$, and the number of such integers is
$$
\sum\limits_{j=1}^k\#\cM_j\ge\sum\limits_{j=1}^k
\frac{0.09N}{q_j}\ge 0.09\eps N.
$$
This completes the proof of Theorem~\ref{thm:PosPropNonres}.

\subsection{Proof of Theorem~\ref{thm:nonres}}

Using Theorem~\ref{thm:PosPropNonres}, we immediately obtain the
following result, which is needed in our proof of
Theorem~\ref{thm:nonres} below:

\begin{lemma}
\label{lem:PosPropXY} Let $\sigma\in\{\pm 1\}$ be fixed. For every
$\eps>0$ there exists a constant $\eta>0$ such that, for all
sufficiently large primes $p$, the lower bound
$$
\#\left\{(n,m)~:~1\le n\le N,~1\le m\le
M,~(nm|p)=\sigma\right\}\ge\eta\,NM
$$
holds with $N=\fl{p^{1/(4\sqrt{e}\,)+\eps}}$ and an arbitrary
positive integer $M$.
\end{lemma}

The next elementary result characterizes the set of values taken by
the Beatty sequence $\cB_{\alpha,\beta}$ in the case that
$\alpha>1$:
\begin{lemma}
\label{lem:Beatty values} Let $\alpha > 1$.  A positive  integer
$m>\beta$ belongs to the Beatty sequence $\cB_{\alpha,\beta}$ if and
only if
$$
0<\{\alpha^{-1}(m-\beta+1)\}\le\alpha^{-1},
$$
and in this case $m = \fl{\alpha n + \beta}$ if and only if $n =
\rf{\alpha^{-1}(m - \beta)}$.
\end{lemma}

The following estimate is a particular case of a series of similar
estimates dating back to the early works of Vinogradov (see, for
example, \cite{Vin2}):

\begin{lemma}
\label{lem:bilinear} Let $\lambda$ be a real number and suppose that
the inequality
$$
\left|\lambda - \frac{r}{q} \right| \le \frac{1}{q^2}
$$
holds for some integers $r$ and $q\ge 1$ with $\gcd(r,q)=1$. Then,
for any complex numbers $a_n,b_m$ such that
$$
\max_{n\le N}\{|a_n|\}\le 1\mand \max_{m\le M}\{|b_m|\}\le 1,
$$
the following bound holds:
$$
\sum_{n\le N}\sum_{m\le M}a_nb_m\,\e(\lambda nm) \ll XY
\sqrt{\frac{1}{X} +\frac{1}{Y} + \frac{1}{q}  + \frac{q}{XY}},
$$
where $\e(z)=\exp(2\pi i z)$ for all $z\in\R$.
\end{lemma}

Considering for every integer $h\ge 1$ the sequence of convergents
in the continued fraction expansion of $\lambda h$, from
Lemma~\ref{lem:bilinear} we derive the following statement:

\begin{corollary}
\label{cor:explicit bilinear} For every irrational $\lambda$, there
are functions $H_\lambda(K)\to\infty$ and $\rho_\lambda(K) \to 0$ as
$K\to\infty$ such that for any complex numbers $a_n,b_m$ such that
$$
\max_{n\le N}\{|a_n|\}\le 1\mand \max_{m\le M}\{|b_m|\}\le 1,
$$
the bound
$$
\left|\sum_{n\le N}\sum_{m\le M}a_nb_m\,\e(\lambda h nm)\right|\le
\rho_\lambda(K) NM
$$
for all integers $h$ in the range $1\le |h|\le H_\lambda(K)$, where
$K=\min\{N,M\}$.
\end{corollary}

In particular, if $\lambda$ is irrational and $h\ne 0$ is fixed,
then
$$
\sum_{n\le N}\sum_{m\le M}a_nb_m\,\e(\lambda h nm)  = o(NM)
$$
whenever $\min\{N,M\}\to\infty$.

\bigskip

We now turn to the proof of Theorem~\ref{thm:nonres}.

\medskip

\textbf{Case 1:} $\alpha>1$. Put $\lambda=\alpha^{-1}$, and let
$\sigma\in\{\pm 1\}$ be fixed. For all integers $N,M\ge 1$ and
primes~$p$, we consider the set of ordered pairs
$$
\cW_p^\sigma(N,M)=\left\{(n,m)~:~1\le n\le N,~1\le m\le
M,~(nm|p)=\sigma\right\}.
$$
For every $\eps>0$, Lemma~\ref{lem:PosPropXY} shows that there is a
constant $\eta>0$ such that, for all sufficiently large primes $p$,
the inequality
$$
\#\cW_p^\sigma(N,M)\ge\eta\,NM
$$
holds with $N=\fl{p^{1/(4\sqrt{e}\,)+\eps/2}}$ and an arbitrary
positive integer $M$. For every large prime $p$, let $N$ be such an
integer, and put $M=\fl{p^{\eps/2}}$. To prove
Theorem~\ref{thm:nonres} when $\alpha>1$, by Lemma~\ref{lem:Beatty
values} it suffices to show that the set
$$
\cV_p^\sigma(N,M)=\left\{(n,m)\in\cW_p^\sigma(N,M)~:~0<\{\lambda
nm-\lambda\beta+\lambda\}\le\lambda\right\}
$$
is nonempty for $\sigma=-1$ when $p$ is sufficiently large. In fact,
we shall prove this result for either choice of $\sigma\in\{\pm
1\}$.

To simplify the notation, write $\cW^\sigma=\cW_p^\sigma(N,M)$ and
$\cV^\sigma=\cV_p^\sigma(N,M)$. To estimate $\#\cV^\sigma$, we use
the well known \emph{Erd\H{o}s--Tur{\'a}n inequality} between the
discrepancy of a sequence and its associated exponential sums; for
example, see~\cite[Theorem~2.5, Chapter~2]{KuNi}. For any integer
$H>1$, we have
$$
\left|\#\cV^\sigma- \lambda\,\#\cW^\sigma\right|\ll \frac{
\#\cW^\sigma}{H}+ \sum_{h=1}^{H}\frac{1}{h}\left|\sum_{(n,m) \in
\cW^\sigma}\e(\lambda h nm)\right|.
$$
Applying Corollary~\ref{cor:explicit bilinear} with the choice
$$
H=\min\left\{H_\lambda(K), \exp\(\rho_\lambda(K)^{-1/2}\)\right\}
$$
where $K=\min\{N,M\}$ as before, we see that
$$
\left|\#\cV^\sigma- \lambda\,\#\cW^\sigma\right|\ll
\frac{\#\cW^\sigma}{H}+  \rho_\lambda(K)\,NM \log H\ll\frac{NM}{\log
H}\,.
$$
Since $H\to\infty$ as $p\to\infty$, and the lower bound
$$
\#\cW^\sigma\ge\eta\,NM
$$
holds by Lemma~\ref{lem:PosPropXY}, it follows that
$$
\#\cV^\sigma\ge(\lambda\,\eta+o(1))\,NM\qquad(p\to\infty).
$$
In particular, $\cV^\sigma\ne\varnothing$ for either choice of
$\sigma\in\{\pm 1\}$ once $p$ is sufficiently large.

\medskip

\textbf{Case 2:} $0<\alpha<1$. In this case,
Theorem~\ref{thm:nonres} follows easily from the classical Burgess
bound for the least quadratic non-residue modulo~$p$ since the
sequence $\cB_{\alpha,\beta}$ contains all integers exceeding
$\fl{\alpha+\beta}$.

\medskip

\textbf{Case 3:} $\alpha<0$. We note that the identity
$$
\fl{\alpha n+\beta}=-\fl{-\alpha n-\beta+1}
$$
holds for all $n\ge 1$ with at most $O(1)$ exceptions (since
$\alpha$ is irrational), hence the sequences $\cB_{\alpha,\beta}$
and $-\cB_{-\alpha,-\beta+1}$ are essentially the same.

If $\alpha<-1$, we argue as in Case~1 with $\alpha$ replaced by
$-\alpha>1$ and $\beta$ replaced by $-\beta+1$. Choosing
$\sigma=-(-1|p)$, Theorem~\ref{thm:nonres} then follows from the
fact that $\cV^\sigma\ne\varnothing$ once $p$ is sufficiently large.

Finally, if $-1<\alpha<0$, we note that the sequence
$\cB_{\alpha,\beta}$ contains all integers up to
$\fl{\alpha+\beta}$. Hence, the result follows from the Burgess
bound in the case that $(-1|p)=+1$ and from the ubiquity of
quadratic residues modulo~$p$ in the case that $(-1|p)=-1$.

\subsection{Proof of Theorem~\ref{thm:nonresPiat}}

The following statement is a variant of~\cite[Lemma 4.3]{GrKol} (we
omit the proof, which follows the same lines):

\begin{lemma}
\label{lem:bilinear2} Let $L$ and $M$ be large positive parameters,
and let $(\kappa,\lambda)$ be an exponent pair. Then for any complex
numbers $a_\ell,b_m$ such that
$$
\max_{L/2<\ell\le L}\{|a_\ell|\}\le 1\mand \max_{M/2<m\le
M}\{|b_m|\}\le 1,
$$
the bound
\begin{eqnarray*}
\lefteqn{\left|\sum_{L/2<\ell\le L}\sum_{M/2<m\le M}a_\ell
b_m\e(h\ell^{1/c}m^{1/c})\right|}\\
&& \quad \ll
\(h^{\kappa_0}L^{\kappa_0/c+\lambda_0}M^{1-\kappa_0+\kappa_0/c}
+h^{-1/2}(LM)^{1-1/(2c)}+LM^{1/2}\)\log L
\end{eqnarray*}
holds for any $h\ge 1$, where
$$
\kappa_0 = \frac{\kappa}{2\kappa+2} \qquad\text{and} \qquad
\lambda_0 =\frac{\kappa+\lambda+1}{2\kappa+2}.
$$
\end{lemma}

Turning to the proof of Theorem~\ref{thm:nonresPiat}, let us fix $c$
in the range
$$
1<c<1+\frac{1-\lambda}{2\kappa-\lambda+3}.
$$
If $L$ and $M$ are sufficiently large, and $\ell\in (L/2, L]$, $m\in
(M/2,M]$ are integers such that
$$
1-\frac{1}{2(LM)^{1-1/c}}\le \{\ell^{1/c}m^{1/c}\},
$$
then $\fl{n^c}=\ell m$ for some integer $n$. Indeed, with
$n=\fl{\ell^{1/c}m^{1/c}}+1$ we see that $n^c\ge\ell m$, and also
\begin{equation*}
\begin{split}
n^c&=\(\ell^{1/c}m^{1/c} +1 - \{\ell^{1/c}m^{1/c}\}\)^c\\
&\le \ell m \( 1 + \frac{1}{2 \ell^{1/c}m^{1/c} (LM)^{1-1/c}}\)^c\\
&\le\ell m \( 1 + \frac{1}{2 \ell m}\)^c<\ell m+1,
\end{split}
\end{equation*}
where the last inequality holds if $L$ and $M$ are large enough.
Below, we work with integers $L,M$ that tend to infinity with the
prime $p$.

Let
$$
J = \rf{\frac{\log (2/\delta)}{\log 2}} \mand
\delta_1=\frac{\delta}{2(J+1)},
$$
where $\delta$ is as in Theorem~\ref{thm:PosPropNonres}. Since
$2^{-J-1}<\delta/2$, by considering the intervals
$(2^{-j-1}p^{1/4\sqrt{e}+\eps}, 2^{-j}p^{1/4\sqrt{e}+\eps}]$ for $j
=0, \ldots, J$ we see that there is an integer $L$ with
$2^{-J}p^{1/4\sqrt{e}+\eps} <L \le p^{1/4\sqrt{e}+\eps}$ such that
the interval $(L/2, L]$ contains a set $\cL$ with $\# \cL \ge
\delta_1 L$ quadratic non-residues modulo $p$. Let $A$ be a large
positive constant. From the aforementioned result of
Hall~\cite{Hall} we see that there exists an integer $M$ with
$$
L^{2(c-1)/(2-c)}(\log L)^A\ll M\ll L^{2(c-1)/(2-c)}(\log L)^A
$$
such that the interval $(M/2, M]$ contains a set $\cM$ with
$\#\cM\ge\delta_2 M$ quadratic residues modulo $p$, where $\delta_2>
0$ is an absolute constant. It suffices to show that for some
integers $\ell\in \cL$, $m\in \cM$ the inequality
$$
1-\frac{1}{2(LM)^{1-1/c}}\le \{\ell^{1/c}m^{1/c}\}
$$
holds. As in the proof of Theorem~\ref{thm:nonres}, from the
Erd\H{o}s--Tur{\'a}n inequality we see that  for any $H \ge 1$ the
number of solutions $T$ of this inequality is
\begin{eqnarray*}
T &  = & \frac{\# \cL \# \cM}{2(LM)^{1-1/c}} +  O\(\frac{LM}{H} +
\sum_{h=1}^{H}\frac{1}{h}
\left|\sum_{L/2<\ell\le L}\sum_{M/2<m\le M}\e(h\ell^{1/c}m^{1/c})\right|\)\\
& \ge & 0.5\delta_1\delta_2(LM)^{1/c}
-\frac{c_0LM}{H}-c_0\sum_{h=1}^{H}\frac{1}{h}
\left|\sum_{L/2<\ell\le L}\sum_{M/2<m\le
M}\e(h\ell^{1/c}m^{1/c})\right|,
\end{eqnarray*}
where $c_0$ is an absolute constant. Take
$H=\rf{4c_0(LM)^{1-1/c}/(\delta_1\delta_2)}$. With this choice it
suffices to prove that
$$
c_0\sum_{h=1}^{H}\frac{1}{h}\left|\sum_{L/2<\ell\le L}\sum_{M/2<m\le
M}\e(h\ell^{1/c}m^{1/c})\right|< 0.1 \delta_1\delta_2(LM)^{1/c}.
$$
If $A$ is large enough, this inequality follows from
Lemma~\ref{lem:bilinear2}, which in turn implies that $T > 0$ and
concludes the proof.

\section{Remarks}

We are grateful to the referee who has pointed that some recent work
of Granville and Soundararajan (unpublished) contains the following
result, which yields a stronger form of our
Theorem~\ref{thm:PosPropNonres}:

\begin{theorem}
\label{thm:GSref} Let $x$ be large, and let $f$ be a completely
multiplicative function with $-1\le f(n)\le 1$ for all $n.$ Suppose
that
$$
\sum_{n\le x}f(n)=o(x).
$$
Then for $1/\sqrt e\le \alpha\le 1$ we have
$$
\left|\sum_{n\le x^{\alpha}}f(n)\right|\le (\max\{|\xi|, 1/2
+2(\log\alpha)^2\}+o(1))x^{\alpha}
$$
where
$$
\xi = 1- 2 \log(1 + \sqrt{e}\,)+4 \int_1^{\sqrt{e}} \frac{\log
t}{t+1}\,dt = -0.656999\cdots \,.
$$
\end{theorem}
We note that $\xi$ is the same constant that appears
in~\cite[Theorem~1]{GrSo} (where it is called $\delta_1$, which has
a different meaning in our paper).

The referee has suggested that the following conjecture seems
natural:

\begin{conjecture}
\label{con:GSref} Let $x$ be large, let $f$ be a completely
multiplicative function with $-1\le f(n)\le 1$ for all $n$. and
suppose that
$$
\sum_{n\le x}f(n)=o(x).
$$
Then for $1/\sqrt e\le \alpha\le 1$ we have
$$
\left|\sum_{n\le x^{\alpha}}f(n)\right|\le
(-2\log\alpha+o(1))x^{\alpha}.
$$
\end{conjecture}

Finally, the referee also observes that Theorem~\ref{thm:nonres}
holds also for rational $\alpha\not=0$. The proof uses recent work
of Balog, Granville and Soundararajan~\cite{BaGrSo}.


\begin{thebibliography}{99}

\bibitem{Abe}
A.~G.~Abercrombie, `Beatty sequences and multiplicative number
theory', \emph{Acta Arith.} \textbf{70} (1995),  195--207.

\bibitem{BaGrSo} A. Balog, A.~Granville and K.~Soundararajan,
`Multiplicative functions in arithmetic progressions',
\emph{Preprint}, 2007 (available from {\tt
http://arxiv.org/abs/math/0702389}).

\bibitem{BaSh1}
W.~Banks  and I.~E.~Shparlinski, `Non-residues and primitive roots
in Beatty sequences', \emph{Bull.\ Austral.\ Math.\ Soc.}
\textbf{73} (2006), 433--443.

\bibitem{BaSh2}
W.~Banks  and I.~E.~Shparlinski, `Short character sums with Beatty
sequences', \emph{Math.\ Res.\ Lett.} \textbf{13} (2006), 539--547.

\bibitem{Beg}
A.~V.~Begunts, `An analogue of the Dirichlet divisor problem',
\emph{Moscow Univ.\ Math.\ Bull.} \textbf{59} (2004),   no.~6,
37--41.

\bibitem{Burg}
D.~A.~Burgess, `The distribution of quadratic residues and
non-residues', \emph{Mathematika} \textbf{4} (1957), 106--112.

\bibitem{Gar}
M.~Z.~Garaev, `A note on the least quadratic non-residue of the
integer-sequences', \emph{Bull.\ Austral.\ Math.\ Soc.} \textbf{68}
(2003), 1--11.

\bibitem{Gr}
S. W. Graham, `An algorithm for computing optimal exponent pairs',
\emph{J. London Math. Soc. (2)} {\textbf 33} (1986), 203--218.

\bibitem{GrKol}  S. W. Graham and  G. Kolesnik,
\emph{Van der Corput's Method of Exponential Sums}, Cambridge Univ.
Press, 1991.

\bibitem{GrRing}  S. W. Graham and  C. J. Ringrose,
`Lower bounds for least quadratic nonresidues', \emph{Analytic
number theory (Allerton Park, IL, 1989)}, Birkh\"auser, Boston, MA,
1990,  269--309.

\bibitem{GrSo} A.~Granville and K.~Soundararajan,
`The spectrum of multiplicative functions', \emph{Ann.\ Math.}
\textbf{153} (2001),   407--470.

\bibitem{Hall} R.~R.~Hall, `Proof of a conjecture of Heath-Brown concerning
quadratic residues', \emph{Proc.\ Edinburgh Math.\ Soc.\ (2)}
\textbf{39} (1996), 581--588.

\bibitem{HarWr} G.~H.~Hardy and E.~M.~Wright,  \emph{An Introduction
to the Theory of Numbers}, 5th ed., Oxford, 1979.

\bibitem{Hild} A.~Hildebrand, `A note on Burgess' character sum
estimate',  \emph{C.\ R.\ Math.\ Rep.\ Acad.\ Sci.\ Canada}
\textbf{8} (1986), no.~1, 35--37.

\bibitem{Huxley}  M. N. Huxley,
\emph{Area, lattice points and exponential sums}, Oxford Univ.
Press, 1996.

\bibitem{Ivic}  I. A. Ivi{\'c},
\emph{The Riemann zeta-function}, John Willey, 1985.

\bibitem{Kom1}
T.~Komatsu, `A certain power series associated with a Beatty
sequence', \emph{Acta Arith.}  \textbf{76}  (1996), 109--129.

\bibitem{Kom2}
T.~Komatsu, `The fractional part of $n\vartheta+\varphi$ and Beatty
sequences',  \emph{J.\ Th\'eor.\ Nombres Bordeaux} \textbf{7}
(1995),  387--406.

\bibitem{KuNi}
L.~Kuipers and H.~Niederreiter, \emph{Uniform distribution of
sequences}, Wiley-Interscience, New York-London-Sydney, 1974.

\bibitem{LW} Y. K. Lau and J. Wu, `On the least quadratic
non-residue', \emph{Preprint}, 2006, {\tt
http://hal.archives-ouvertes.fr/hal-00097136/en/}.

\bibitem{OB}
K.~O'Bryant, `A generating function technique for Beatty sequences
and other step sequences', \emph{J.\ Number Theory} \textbf{94}
(2002),  299--319.

\bibitem{Preo1}
S.~N.~Preobrazhenski\u\i, `On the least quadratic non-residue in an
arithmetic sequence', \emph{Moscow Univ.\ Math.\ Bull.} \textbf{56}
(2001), no.~1, 44--46.

\bibitem{Preo2}
S.~N.~Preobrazhenski\u\i, `On power non-residues modulo a prime
number in a special integer sequence', \emph{Moscow Univ.\ Math.\
Bull.} \textbf{56} (2001), no.~4, 41--42.

\bibitem{Preo3}
S.~N.~Preobrazhenski\u\i, `On the least power non-residue in an
integer sequence', \emph{Moscow Univ.\ Math.\ Bull.} \textbf{59}
(2004), no.~1, 33--35.

\bibitem{Rob1} O. Robert and P. Sargos, `A fourth
derivative test for exponential sums', \emph{Compositio Math.} {\bf
130} (2002), no. 3, 275--292.

\bibitem{Rob2} O. Robert and P. Sargos, `A third
derivative test for mean values of exponential sums with application
to lattice point problems',  \emph{Acta Arith.} {\bf 106} (2003),
no. 1, 27--39.

\bibitem{Sarg} P. Sargos,  `An analog of van der
Corput's $A\sp 4$-process for exponential sums',  \emph{Acta Arith.}
{\bf 110} (2003), no. 3, 219--231.

\bibitem{Tijd}
R.~Tijdeman, `Exact covers of balanced sequences and Fraenkel's
conjecture', \emph{Algebraic number theory and Diophantine analysis
(Graz, 1998)}, 467--483, de Gruyter, Berlin, 2000.

\bibitem{Vin2} I. M. Vinogradov,
`An improvement of the estimation of sums with primes',
\emph{Izvestia Akad. Nauk SSSR} \textbf{17} (1943), 17--34 (in
Russian).

\end{thebibliography}
\end{document}